\documentclass{article}

\usepackage{amssymb}
\usepackage{amsmath}
\usepackage[dvips]{graphicx}
\begin{document}

%%% remove comment delimiter ('%') and select language if required
%\selectlanguage{spanish}

\textbf{$ $  The Cauchy Problem Of  The Moment Theoiry Elasticity In $R^{n} $}

I. E. Niyozov, O. I. Makhmudov

 \textbf{ Abstract}

In this paper, we considered the problem of analytical continuation of the solution of the system equations of the moment theory of elasticity in spacious bounded domain from its values and values of its strains on part of the boundary of this domain, i.e., the Cauchy's problem.

\textbf{ Key words}: the Cauchy problem, system theory of elasticity, elliptic system, ill-posed problem, Carleman matrix, regularization.

 \textbf{ 1. Introduction}

In this paper, we considered the problem of analytical continuation of the solution of the system equations of the moment theory of elasticity in spacious bounded domain from its values and values of its strains on part of the boundary of this domain, i.e., the Cauchy's problem.

Since, in many actual problems, either a part of the boundary is inaccessible for measurement of displacement and tensions or only some integral characteristic are available. In experimental study of the stress-strain state of actual constructions, we can make measurements only on the accessible part of the surface.

In a practical investigation of experimental dates or diagnostic moving abject arise problems of estimation concerning deformed position of the object. Solution of the problems by using well known classical propositions is connected to difficulties of absence of experimental dates which is necessary for formulation of boundary value (classical) conditions.

Therefore it is necessary consider the problem of continuation for solution of elasticity system of equations to the domain by values of solutions and normal derivatives in the part of boundary of domain.

System equation of moment theory elasticity is elliptic. Therefore the problem Cauchy for this system is ill-posed. For ill-posed problems, one does not prove the existence theorem: the existence is assumed a priori. Moreover, the solution is assumed to belong to some given subset of the function space, usually a compact one [1]. The uniqueness of the solution follows from the general Holmgren theorem [2]. On establishing uniqueness in the article studio of ill-posed problems, one comes across important questions concerning the derivation of estimates of conditional stability and the construction of regularizing operators.

Our aim is to construct an approximate solution using the Carleman function method.

Let $x=(x_{1} ,...,x_{n} )$ and $y=(y_{1} ,...,y_{n} )$ be points of the $n$-dimensional Euclidean space $E^{n} $, $D$ a bounded simply connected domain in $E^{n} $, with piecewise-smooth boundary consisting of a piece $\Sigma $ of the plane $y_{n} =0$ and a smooth surface $S$ lying in the half-space $y_{n} >0.$

Suppose that $2n$-component vector function

$U(x)=(u_{1} (x),...,u_{n} (x),w_{1} (x),...,w_{n} (x))=(u(x),w(x))$ satisfied in $D$ the system equations moments theory elasticity [3]:

\[\left\{\begin{array}{llll} {(\mu +\alpha )\Delta u+(\lambda +\mu -\alpha )graddivu+2\alpha \, rotw+\rho \sigma ^{2} u=0,} & {} & {} & {} \\ {(\nu +\beta )\Delta w+(\varepsilon +\nu -\beta )graddivw+2\alpha \, rotu-4\alpha w+\theta \sigma ^{2} w=0,} & {} & {} & {} \end{array}\right. (1),\]
where $\lambda ,\, {\kern 1pt} \mu ,\, {\kern 1pt} \nu ,\, {\kern 1pt} \beta ,{\kern 1pt} \, \varepsilon ,{\kern 1pt} \, \alpha ,\rho ,\, \sigma $ is coefficients which characterizing medium, satisfying the conditions

\[\mu >0,{\kern 1pt} \, \, 3\lambda +2\mu >0,{\kern 1pt} \, \, \alpha >0,{\kern 1pt} \, \, \varepsilon >0,{\kern 1pt} \, \, 3\varepsilon +2\nu >0,\, {\kern 1pt} \, \beta >0,\, \rho >0,\, \sigma >0.\]

For brevity it is convenient to use matrix notation. Let us introduce the matrix differential operator

\[M(\partial _{x} )=\left|\left|\begin{array}{ccc} {M^{(1)} } & {M^{(2)} } & {} \\ {M^{(3)} } & {M^{(4)} } & {} \end{array}\right|\right|,\]
where

\[M^{(i)} =\left\| M_{k\, j}^{(i)} \right\| _{n\times n} ,\quad i=1,2,3,4,\]
moreover

\[M_{k\, j}^{(1)} =\delta _{k\, j} (\mu +\alpha )(\Delta +\sigma _{1}^{2} )+(\lambda +\mu -\alpha )\frac{\partial ^{2} }{\partial x_{k} \partial x_{j} } ,\, \, \, \, \, k,j=1,...,n\]

\[M_{k\, j}^{(2)} =M_{k\, j}^{(3)} =-2\alpha \sum _{p=1}^{n} \varepsilon _{k\, j\, p} \frac{\partial }{\partial x_{p} } ,\, \, \, k,j=1,...,n,\]

\[M_{k\, j}^{(4)} =\delta _{k\, j} \left[(\nu +\beta )\Delta +\sigma _{2}^{2} \right]+(\varepsilon +\nu -\beta )\frac{\partial ^{2} }{\partial x_{k} \partial x_{j} } ,\, \, \, k,j=1,...,n,\]
here

\[\sigma _{1}^{2} =\frac{\rho \sigma ^{2} }{\mu +\alpha } ,\, \, \, \sigma _{2}^{2} =\frac{\theta \sigma ^{2} -4\alpha }{\nu +\beta }  , \delta _{k\, j} =\left\{\begin{array}{llll} {1,\quad if\quad k=j} & {} & {} & {} \\ {0,\quad if\quad k\not =j,} & {} & {} & {} \end{array}\right. \]
$\varepsilon _{k\, j\, p} $ so-called $\varepsilon -$tensor or Levi-Civita's symbol, which defend following equaliti's

\[\varepsilon _{k\, j\, p} =\left\{\begin{array}{llll} {0,\quad if\, at\, \, least\, \, two\, \, of\, \, three-subscripts\quad k,j,p\quad are\, \, equal,} & {} & {} & {} \\ {1,\quad if\quad (k,j,p)\quad is\, \, an\, \, even\, permutation,} & {} & {} & {} \\ {-1,\quad if\quad (k,j,p)\quad is\, an\, \, odd\, \, permutation.} & {} & {} & {} \end{array}\right. \]
Then system \eqref{GrindEQ__1_} maybe write in matrix from in the following way:

\[M(\partial _{x} )U(x)=0\, \, \, \, \, \, \, \, \, \, \, \, \, \, \, \, \, \, \, \, \, \, \, \, \, \, \, \, \, \, \, \, \, \, \, \, \, \, \, \, \, \, \, \, \, \, \, \, \, \, \, \, \, (2)\]

A solution $U$ of system \eqref{GrindEQ__1_} in the domain $D$ is said to be regular if $U\in C^{1} (\overline{D})\bigcap _{} C^{2} (D).$

\textbf{ Statement of the problem.} Find a regular solution $U$ of system \eqref{GrindEQ__1_} in the domain $D$ using its Cauchy data on the surface $S$:

\[U(y)=f(y),\quad \quad T(\partial _{y} ,n(y))U(y)=g(y),\quad y\in S,\, \, \, \, \, \, \, \, \, \, \, \, \, \, \, \, \, \, \, \, \, (3)\]
where $T(\partial _{y} ,n(y))$ is the stress operator, i.e.,

\[T(\partial _{y} ,n(y))=\left|\left|\begin{array}{ccc} {T^{(1)} (\partial _{y} ,n)} & {T^{(2)} (\partial _{y} ,n)} & {} \\ {T^{(3)} (\partial _{y} ,n)} & {T^{(4)} (\partial _{y} ,n)} & {} \end{array}\right|\right|,\]

\[T^{(i)} (\partial _{y} ,n)=\left\| T_{k\, j}^{(i)} (\partial _{y} ,n)\right\| _{n\times n} ,\quad i=1,2,3,4,\]

\[T_{k\, j\, }^{(1)} (\partial _{y} ,n)=\lambda \, n_{k} \frac{\partial }{\partial y_{j} } +(\mu -\alpha )\, n_{j} (y)\frac{\partial }{\partial y_{k} } +(\mu +\alpha )\, \delta _{k\, j} \frac{\partial }{\partial n(y)} ,\]

\[T_{k\, j}^{(2)} (\partial _{y} ,n)=2\alpha \sum _{p=1}^{n} \varepsilon _{k\, j\, p} n_{p} (y),\quad T_{k\, j}^{(3)} (\partial _{y} ,n)=0,\]

\[T_{k\, j}^{(4)} (\partial _{y} ,n)=\varepsilon \, n_{k} (y)\frac{\partial }{\partial y_{j} } +(\nu -\beta )\, n_{j} (y)\frac{\partial }{\partial y_{k} } +(\nu +\beta )\frac{\partial }{\partial n(y)} ,\]

$n(y)=(n_{1} (y),...,n_{n} (y))$ is the unit outward normal vector on $\partial D$ at a point $y$, $f=(f_{1} ,\ldots ,f_{2n} ),$ $g=(g_{1} ,\ldots ,g_{2n} )$ are given continuous vector functions on $S.$

 \textbf{2. Construction of the matrix Carleman and approximate solution for the domain type's cap}

It is well known, that any regular solution $U(x)$ system \eqref{GrindEQ__1_} is specified by the formula

\[U(x)=\int _{\partial D} (\Psi (y,x)\{ T(\partial _{y} ,n)U(y)\} -\{ T(\partial _{y} ,n)\Psi (y,x)\} ^{*} U(y))ds_{y} ,{\kern 1pt} \; x\in D,\, \, \, \, \, \, \, (4)\]
where symbol  is denote of operation transposition, $\Psi (y,x)$ matrix of fundamental solutions system equation of steady-state oscillations of the couple-stress theory of elasticity:

\noindent

\[\Psi (y,x)=\left|\left|\begin{array}{ccc} {\Psi ^{(1)} (y,x)} & {\Psi ^{(2)} (y,x)} & {} \\ {\Psi ^{(3)} (y,x)} & {\Psi ^{(4)} (y,x)} & {} \end{array}\right|\right|,\]
where

\[\Psi ^{(i)} (y,x)=\left\| \Psi _{k\, j}^{(i)} (y,x)\right\| _{n\times n} ,\quad i=1,2,3,4,\]

\[\Psi _{k\, j}^{(1)} (y,x)=\sum _{l=1}^{4}(\delta _{k\, j} \alpha _{l}  +\beta _{l} \frac{\partial ^{2} }{\partial \, x_{k} \partial \, x_{j} } )\varphi _{n} (ik_{l} r),\, \, \, \, \, \, k,j=1,...,n,\]

\[\Psi _{k\, j}^{(2)} (y,x)=\Psi _{k\, j}^{(3)} (y,x)=\frac{2\alpha }{\mu +\alpha } \sum _{l=1}^{4} \sum _{p=1}^{n} \varepsilon _{l} \varepsilon _{k\, j\, p} \frac{\partial }{\partial x_{p} } \varphi _{n} (ik_{l} r),\, \, k,j=1,...,n,\]

\[\Psi _{k\, j}^{(4)} (y,x)=\sum _{l=1}^{4}(\delta _{k\, j} \gamma _{l}  +\delta _{l} \frac{\partial ^{2} }{\partial \, x_{k} \partial \, x_{j} } )\varphi _{n} (ik_{l} r)\, ,\, \, \, \, \, \, k,j=1,...,n,\]
here $\varphi _{n} $ -fudamental solution Helmholtz equation, $r=\left|x-y\right|$,

\[\alpha _{l} =\frac{(-1)^{l} (\sigma _{2}^{2} -k_{l}^{2} )(\delta _{3\, l} +\delta _{4\, l} )}{2\pi (\mu +\alpha )(k_{3}^{2} -k_{4}^{2} )} ,\, \, \, \, \beta _{l} =-\frac{\delta _{1\, l} }{\begin{array}{l} {2\pi \rho \sigma ^{2} } \\ {} \end{array}} +\frac{\alpha _{l} }{k_{l}^{2} } ,\, \, \, \sum _{l=1}^{4}\beta _{l}  =0\]

\[\gamma _{l} =\frac{(-1)^{l} (\sigma _{1}^{2} -k_{l}^{2} )(\delta _{3\, l} +\delta _{4\, l} )}{2\pi (\beta +\nu )(k_{3}^{2} -k_{4}^{2} )} ,\, \, \, \, \delta _{l} =-\frac{\delta _{2\, l} }{\begin{array}{l} {2\pi (\theta \sigma ^{2} -4\alpha )} \\ {} \end{array}} +\frac{\gamma _{l} }{k_{l}^{2} } ,\, \, \, \, \, \, \, \sum _{l=1}^{4}\beta _{l}  =0,\]

\[\varepsilon _{l} =\frac{(-1)^{l} (\delta _{3\, l} +\delta _{4\, l} )}{2\pi (\beta +\nu )(k_{3}^{2} -k_{4}^{2} )} ,\, \, \, \, \sum _{l=1}^{4}\varepsilon _{l}  =0,\, \, \, k_{1}^{2} =\frac{\rho \sigma ^{2} }{\lambda +2\mu } ,\, \, \, k_{2}^{2} =\frac{\theta \sigma ^{2} -4\alpha }{\varepsilon +2\nu } ,\, \,\]
$$
 k_{1}^{2} +k_{2}^{2} =\sigma _{1}^{2} +\sigma _{2}^{2} +\frac{4\alpha ^{2} }{(\mu +\alpha )(\beta +\nu )} ,  k_{1}^{2} k_{2}^{2} =\sigma _{1}^{2} \sigma _{2}^{2} .
$$
Easily we can verity, that $u=\Psi _{j}^{\eqref{GrindEQ__1_}}
(y,x),$ $w=\Psi _{j}^{\eqref{GrindEQ__3_}} (y,x)$ or $u=\Psi
_{j}^{\eqref{GrindEQ__2_}} (y,x),$ $w=\Psi
_{j}^{\eqref{GrindEQ__4_}} (y,x)$ are solution system
\eqref{GrindEQ__1_}, where $\Psi _{j}^{i} (y,x)-$ $j-$ vector
tuple $i-$ matrix.

\textbf{Definition.} \textit{ By the Carleman matrix of problem \eqref{GrindEQ__1_},\eqref{GrindEQ__3_} we mean an $2n\times 2n$ matrix $\Pi (y,x,\tau )$ depending on the two points $y,x$ and positive numerical number parameter $\tau $ satisfying the following two conditions: }

\[1){\kern 1pt} \; \Pi (y,x,\tau )=\Psi (y,x)+G(y,x,\tau ),\]
where matrix $G(y,x,\tau )$ satisfies system \eqref{GrindEQ__1_} with respect to the variable $y$ in the domain $D$, and $\Psi (y,x)$ is a matrix of the fundamental solutions of system \eqref{GrindEQ__1_};

\[2){\kern 1pt} \; \int _{\partial D\backslash S} \left(|\Pi (y,x,\tau )|+|T(\partial _{y} ,n)\Pi (y,x,\tau )|\right)ds_{y} \le \varepsilon (\tau ),\]

where $\varepsilon (\tau )\to 0,$ as $\tau \to \infty ;$ here $|\Pi |$ is the Euclidean norm of the matrix $\Pi =||\Pi _{i\, j} ||_{2n\times 2n} ,$ i.e., $|\Pi |=(\sum _{i,j=1}^{2n} \Pi _{i\, j}^{2} )^{\frac{1}{2} } .$ In particular, $|U|=\left(\sum _{m=1}^{n} (u_{m}^{2} +w_{m}^{2} )\right)^{\frac{1}{2} } .$

It is well known, that for the regular vector functions $v(y)$ and $u(y)$ holds formula [4]:

\[\int _{D} [v(y)\{ M(\partial _{y} )u(y)\} -u(y)\{ M(\partial _{y} )v(y)\} dy=\]

\[=\int _{\partial D} [v(y)\{ T(\partial _{y} ,n)u(y)\} -u(y)\{ T(\partial _{y} ,n)v(y)\} ]ds_{y} .\]
Substituting in this equality $v(y)=G(y,x,\tau )$ and $u(y)=U(y)$ is solution system \eqref{GrindEQ__1_}, we have

\[0=\int _{\partial D} [G(y,x,\tau )\{ T(\partial _{y} ,n)U(y)\} -\{ T(\partial _{y} ,n)G(y,x,\tau )\} ^{*} U(y)]ds_{y} .(5)\]

Now adding \eqref{GrindEQ__4_} and \eqref{GrindEQ__5_}, we have

\textbf{Theorem 1.} \textit{ Any regular solution $U(x)$ of system \eqref{GrindEQ__1_} in the domain $D$ is specified by the formula }

\[U(x)=\int _{\partial D} (\Pi (y,x,\tau )\{ T(\partial _{y} ,n)U(y)\} -\{ T(\partial _{y} ,n)\Pi (y,x,\tau )\} ^{*} U(y))ds_{y} ,{\kern 1pt} \; x\in D,(6)\]
where $\Pi (y,x,\tau )$ is matrix Carleman.

Using the matrix Carleman, easily conclude the estimate stability of solution of the problem \eqref{GrindEQ__1_}, \eqref{GrindEQ__3_} and also indicate effective method decision this problem.

With a view to construct approximate solution of the problem \eqref{GrindEQ__1_}, \eqref{GrindEQ__3_} we construct the following matrix:

\[\Pi (y,x)=\left|\left|\begin{array}{ccc} {\Pi ^{(1)} (y,x)} & {\Pi ^{(2)} (y,x)} & {} \\ {\Pi ^{(3)} (y,x)} & {\Pi ^{(4)} (y,x)} & {} \end{array}\right|\right|,\, \, \, \, \, \, \, \, \, \, \, \, \, \, \, \, \, \, \, \, \, \, \, \, \, \, \, (7)\]

\[\Pi ^{(i)} (y,x)=\left\| \Pi _{k\, j}^{(i)} (y,x)\right\| _{n\times n} ,\quad i=1,2,3,4,\]

\[\Pi _{k\, j}^{(1)} (y,x)=\sum _{l=1}^{4}(\delta _{k\, j} \alpha _{l}  +\beta _{l} \frac{\partial ^{2} }{\partial \, x_{k} \partial \, x_{j} } )\Phi (y,x,k_{l} )\, ,\, \, \, \, \, \, k,j=1,...,n\]

\[\begin{array}{l} {\, \, \, \, \, \, \, \, \, \, \, \, \, \, \, \, \, \, \, \, \, \, \, \, \, \, \, \, \, \, \, \, \, \, \, \, \, \, \Pi _{k\, j}^{(2)} (y,x)=\Pi _{k\, j}^{(3)} (y,x)=} \\ {} \\ {\, \, \, \, \, \, \, \, \, \, \, \, \, \, \, \, \, =\frac{2\alpha }{\mu +\alpha } \sum _{l=1}^{4} \sum _{p=1}^{n} \varepsilon _{l} \varepsilon _{k\, j\, p} \frac{\partial }{\partial x_{p} } \Phi (y,x,k_{l} )\, ,\, \, \, \, \, \, k,j=1,...,n,} \end{array}\]

\[\Pi _{k\, j}^{(4)} (y,x)==\sum _{l=1}^{4}(\delta _{k\, j} \gamma _{l}  +\delta _{l} \frac{\partial ^{2} }{\partial \, x_{k} \partial \, x_{j} } )\Phi (y,x,k_{l} )\, ,\, \, \, \, \, \, \, \, k,j=1,...,n\, \, \, \, \, \, \, \, \, \, \, \, \, \, \, (8)\]

\noindent where

\[C_{n} K(x_{_{n} } )\Phi (y,x,k)=\int _{0}^{\infty } Im[\frac{K(i\sqrt{u^{2} +s} +y_{n} )}{i\sqrt{u^{2} +s} +y_{n} -x_{n} } ]\frac{\psi (ku)\, du}{\sqrt{u^{2} +s} } ,\, \, \, \, \, \, \, \, \, \, \, \, \, \, \, \, \, \, (9)\]

$\psi (ku)=\left\{\begin{array}{l} {uJ_{0} (ku),\, \, \, n=2m,\, \, m\ge 1,} \\ {\cos ku,\, \, \, \, \, \, \, n=2m+1,\, \, \, m\ge 1,} \end{array}\right. $ $J_{0} (u)$-Bessel function of order zero,

\[s=(y_{1} -x_{1} )^{2} +...+(y_{n-1} -x_{n-1} )^{2} ,         C_{2} =2\pi \]

\[C_{n} =\left\{\begin{array}{l} {(-1)^{m} \cdot 2^{-n} (n-2)\pi \omega _{n} (m-2)!,\, \, \, \, \, n=2m} \\ {(-1)^{m} \cdot 2^{-n} (n-2)\pi \omega _{n} (m-1)!,\, \, \, \, \, n=2m+1.} \end{array}\right. \]

            $K(\omega ),{\kern 1pt} \; \omega =u+iv$ ($u,{\kern 1pt} \; v$ are real), is an entire function taking real values on the real axis and satisfying the conditions $K(u)\ne \infty ,\, \, \, \, \, \left|u\right|<\infty ,$

$K(u)\not =0,{\kern 1pt} \; \mathop{\sup }\limits_{v\ge 1} |\exp \nu \, \left|Imk\right|K^{(p)} (\omega )|=M(p,u)<\infty ,{\kern 1pt} \; p=0,...,m,{\kern 1pt} \; u\in R^{1} .$In work [4] proved.

\textbf{Lemma 1.} \textit{ For  function $\Phi (y,x,k)$ the  formula is valid }

\[{\kern 1pt} {\kern 1pt} C_{n} \Phi (y,x,k)=\varphi _{n} (ikr)+g_{n} (y,x,k),\quad r=|y-x|,\, \, \, \, \, \, \, \, \, \, \, \, \, \, \, \, \, \, \, \, \, \, \, \, \, \, \, \, \, \, \, \, \, \, \, \, \, \, (10)\]
where $\varphi _{n} $ -fudamental solution Helmholtz equation, $g_{n} (y,x,k)$ is a regular function that is defined for all $y$ and $x$ satisfies Helmholtz equation: $\Delta (\partial _{y} )g_{n} -k^{2} g_{n} =0.$

\noindent            In \eqref{GrindEQ__9_} we assume the function $K(\omega )=\exp (\tau \omega )$. Then

\[\Phi (y,x,k)=\Phi _{\tau } (y-x,k),  \]

\[C_{n} \Phi _{\tau } (y-x,k)=\frac{\partial ^{m-1} }{\partial s^{m-1} } \int _{0}^{\infty } Im[\frac{\exp \, \tau (i\sqrt{u^{2} +s} +y_{n} -x_{n} )}{i\sqrt{u^{2} +s} +y_{n} -x_{n} } ]\frac{\psi (ku)\, du}{\sqrt{u^{2} +s} } =\]

\[=\exp \, \tau (y_{n} -x_{n} )\, \frac{\partial ^{m-1} }{\partial s^{m-1} } \int _{0}^{\infty } \left[\right. -\cos \tau \sqrt{u^{2} +\alpha ^{2} } +\]

\begin{equation} \label{GrindEQ__11_}
+(y_{n} -x_{n} )\frac{\sin \tau \sqrt{u^{2} +s} }{\sqrt{u^{2} +s} } \left. \right]\, \psi (ku)du,
\end{equation}

\[\Phi '_{\tau } (y-x,k)=\frac{\partial \Phi _{\tau } }{\partial \tau } .\]

\[C_{n} \Phi '_{\tau } (y-x,k)=\exp \, \tau (y_{n} -x_{n} )\, \frac{\partial ^{m-1} }{\partial s^{m-1} } \int _{0}^{\infty } \frac{\sin \tau \sqrt{u^{2} +s} }{\sqrt{u^{2} +s} } \psi (ku)du,\]

\begin{equation} \label{GrindEQ__12_}
C_{n} \Phi '_{\tau } (y-x,k)=\exp \, \tau (y_{n} -x_{n} )\, \frac{\partial ^{m-1} }{\partial s^{m-1} } \psi '_{\tau } (k,s),
\end{equation}

\[\psi '_{\tau } =\left\{\begin{array}{l} {0,\, \, \, \, \, \, \, \, \, \, \, \, \, \, \, \, \, \, \, \, \, \, \, \, \, \, \, \, \, \, \, \, \, \, \, \, \, \, \, \, \, \, \, \tau <k} \\ {\cos \sqrt{s(\tau ^{2} -k^{2} )} ,\, \, \, \, \, \, \, \, \, \, \, \, \, \, \, n=2m\, \, \, } \\ {\frac{1}{2} \pi \, J{}_{0} (\sqrt{s(\tau ^{2} -k^{2} )} ),\, \, \, \, \, \, \tau >k} \end{array}\right. \]

Now  in formul \eqref{GrindEQ__7_}, \eqref{GrindEQ__8_} and \eqref{GrindEQ__9_} to take  $\Phi (y,x,k)=\Phi _{\tau } (y-x,k)$, we construct matrix  $\Pi (y,x)=\Pi (y,x,\tau )$

From Lemma 1 we obtain.

\textbf{ Lemma 2.} \textit{ The matrix $\Pi (y,x,\tau )$ given by \eqref{GrindEQ__7_} and \eqref{GrindEQ__8_} is Carleman's matrix for problem \eqref{GrindEQ__1_},\eqref{GrindEQ__3_}.}

\textbf{ Proof.} By \eqref{GrindEQ__7_}, \eqref{GrindEQ__8_}, \eqref{GrindEQ__9_} and Lemma 1 we have

\[\Pi (y,x,\tau )=\Psi (y,x)+G(y,x,\tau ),\]
where

\[G(y,x,\tau )=\left\| \begin{array}{ccc} {G^{(1)} (y,x,\tau )} & {G^{(2)} (y,x,\tau )} & {} \\ {G^{(3)} (y,x,\tau )} & {G^{(4)} (y,x,\tau )} & {} \end{array}\right\| ,\]

\[G^{(i)} (y,x,\tau )=\left\| G_{k\, j}^{(i)} (y,x,\tau )\right\| _{n\times n} ,\quad i=1,2,3,4,\]

\[G_{k\, j}^{(1)} (y,x,\tau )=\sum _{l=1}^{4}(\delta _{k\, j} \alpha _{l}  +\beta _{l} \frac{\partial ^{2} }{\partial \, x_{k} \partial \, x_{j} } )g_{n} (y,x,k_{l} ,\tau )\, ,\, \, \, \, \, \, k,j=1,...,n\]

\[\begin{array}{l} {\, \, \, \, \, \, \, \, \, \, \, \, \, \, \, \, \, \, \, \, \, \, \, \, \, \, \, \, \, \, \, \, \, \, \, \, \, \, G_{k\, j}^{(2)} (y,x,\tau )=G_{k\, j}^{(3)} (y,x,\tau )=} \\ {} \\ {\, \, \, \, \, \, \, \, \, \, \, \, \, \, \, \, \, =\frac{2\alpha }{\mu +\alpha } \sum _{l=1}^{4} \sum _{p=1}^{n} \varepsilon _{l} \varepsilon _{k\, j\, p} \frac{\partial }{\partial x_{p} } g_{n} (y,x,k_{l} ,\tau )\, ,\, \, \, \, \, \, k,j=1,...,n,} \end{array}\]

\[G_{k\, j}^{(4)} (y,x,\tau )==\sum _{l=1}^{4}(\delta _{k\, j} \gamma _{l}  +\delta _{l} \frac{\partial ^{2} }{\partial \, x_{k} \partial \, x_{j} } )g_{n} (y,x,k_{l} ,\tau )\, ,\, \, \, \, \, \, \, \, k,j=1,...,n\, \, \, \, \, \, \, \, \, \, \, \, \, \, \]
By a straightforward calculation, we can verify that the matrix $G(y,x,\tau )$ satisfies system \eqref{GrindEQ__1_} with respect to the variable $y$ everywhere in $D$. By using \eqref{GrindEQ__7_}, \eqref{GrindEQ__8_},\eqref{GrindEQ__11_} and \eqref{GrindEQ__12_} we obtain

\[\int _{\partial D\backslash S} \left(|\Pi (y,x,\tau )|+|T(\partial _{y} ,n)\Pi (y,x,\tau )|\right)\, ds_{y} \le C_{1} (x)\, \tau ^{m} \exp (-\tau \, x_{n} ),\, \, \, \, \, \, (13)\]
where $C_{1} (x)$ some bounded function inside $D.$ The lemma is thereby proved.

Let us set

\[U_{\tau } (x)=\int _{S} [\Pi (y,x,\tau )\{ T(\partial _{y} ,n)U(y)\} -\{ T(\partial _{y} ,n)\Pi (y,x,\tau )\} ^{*} U(y)]ds_{y} .\, \, (14)\]

The following theorem holds.

\textbf{Theorem 1.} \textit{ Let $U(x)$ be a regular solution of system \eqref{GrindEQ__1_} in $D$ such that }

\[|U(y)|+|T(\partial _{y} ,n)U(y)|\le M,{\kern 1pt} \; y\in \partial D\backslash S.\, \, \, \, \, \, \, \, \, \, \, \, \, \, \, \, \, \, \, \, \, \, \, \, \, \, \, \, \, \, \, \, \, \, \, \, \, \, \, \, \, \, \, (15)\]
Then for $\tau \ge 1$ the following estimate is valid:

\[|U(y)-U_{\tau } (y)|\le MC_{2} (x)\tau ^{m} \exp (-\tau \, x_{n} ).\]

\textbf{Proof.} By formula \eqref{GrindEQ__6_} and \eqref{GrindEQ__14_}, we have

\[|U(x)-U_{\tau } (x)|=\int _{\partial D\backslash S} [\Pi (y,x,\tau )\{ T(\partial _{y} ,n)U(y)\} -\{ T(\partial _{y} ,n)\Pi (y,x,\tau )\} ^{*} U(y)]ds_{y} .\]
Now on the basis of \eqref{GrindEQ__13_} and \eqref{GrindEQ__15_} we obtain required inequality. The theorem is thereby proved.

Now we write out a result that allows us to calculate $U(x)$ approximately if, instead of $U(y)$ and $T(\partial _{y} ,n)U(y),$ their continuous approximations $f_{\delta } (y)$ and $g_{\delta } (y)$ are given on the surface $S$:

\[\mathop{\max }\limits_{S} |f(y)-f_{\delta } (y)|+\mathop{\max }\limits_{S} |T(\partial _{y} ,n)U(y)-g_{\delta } (y)|\le \delta ,{\kern 1pt} \; 0<\delta <1.\, \, \, \, \, \, \, \, (16)\]
          We define a function $U_{\tau \, \delta } (x)$ by setting

\[U_{\tau \, \delta } (x)=\int _{S} [\Pi (y,x,\tau )g_{\delta } (y)-\{ T(\partial _{y} ,n)\Pi (y,x,\tau )\} ^{*} f_{\delta } (y)]ds_{y} ,\, \, \, \, \, \, \, \, \, \, \, \, \, (17)\]
where

\[\tau =\frac{1}{x_{n}^{0} } ln\frac{M}{\delta } ,\quad x_{n}^{0} =\mathop{\max }\limits_{D} x_{n} ,\quad x_{n} >0.\]

\textbf{Theorem 2.} \textit{ Let $U(x)$ be a regular solution of system \eqref{GrindEQ__1_} in $D$ satisfying condition \eqref{GrindEQ__15_}. Then the following estimate is valid: }

\[|U(x)-U_{\tau \, \delta } (x)|\le C_{3} (x)\delta ^{\frac{x_{n} }{x_{n}^{0} } } \left(ln\frac{M}{\delta } \right)^{m} ,\quad x\in D.\]

\textbf{Proof.} From formula \eqref{GrindEQ__6_} and \eqref{GrindEQ__17_} we have

\[U(x)-U_{\tau \, \delta } (x)=\int _{S} [\Pi (y,x,\tau )\{ T(\partial _{y} ,n)U(y)-g_{\delta } (y)\} -\]

\[-\{ T(\partial _{y} ,n)\Pi (y,x,\tau )\} (U(y)-f_{\delta } (y))]ds_{y} +\]

\[+\int _{\partial D\backslash S} [\Pi (y,x,\tau )\{ T(\partial _{y} ,n)U(y)\} -\{ T(\partial _{y} ,n)\Pi (y,x,\tau )\} ^{*} U(y)]ds_{y} .\]
By the assumption of the theorem and inequalities \eqref{GrindEQ__13_}, \eqref{GrindEQ__15_} and \eqref{GrindEQ__16_} for the any $x\in D,$ we obtain

\[|U(x)-U_{\tau \, \delta } (x)|=C_{2} (x)\delta \tau ^{m} \exp \tau (x_{n}^{0} -x_{n} )+C_{1} (x)\tau ^{m} \exp (-\tau \, x_{n} )\le \]

\[\le C_{3} (x)\tau ^{m} (M+\delta \exp \tau \, x_{n}^{0} )\exp (-\tau \, x_{n} ).\]
Now, it to take $\tau =\frac{1}{x_{n}^{0} } ln\frac{M}{\delta } ,$ then we obtain to proof theorem. The theorem is thereby proved.

\textbf{Theorem 3.} \textit{ Let $U(x)$ be a regular solution of system \eqref{GrindEQ__1_} in $D$ satisfying conditions }

\[|U(y)|+|T(\partial _{y} ,n)U(y)|\le M,{\kern 1pt} \; y\in \partial D\backslash S,\]

\[|U(y)|+|T(\partial _{y} ,n)U(y)|\le \delta ,\quad 0<\delta <1,\quad y\in S.\]
Then

\[|U(x)|\le C_{4} (x)\delta ^{\frac{x_{n} }{x_{n}^{0} } } \left(ln(\frac{M}{\delta } )\right)^{m} ,\]
where $C_{4} (x)=\tilde{C}\int _{\partial D} \frac{1}{r^{n} } ds_{y} ,\quad \tilde{C}-$ constant depending on $\lambda ,\mu ,\varepsilon ,\beta ,\nu .$

\textbf{Proof.} On the basis of Theorem 2 we obtain

\[|U(x)|\le |U_{\tau } (x)|+MC_{2} (x)\tau ^{m} \exp (-\tau \, x_{n} ).\]
Next from the condition theorem and \eqref{GrindEQ__7_}, \eqref{GrindEQ__12_} we obtain

\[|U_{\tau \, \delta } (x)|=\left|\int _{S} [\Pi (y,x,\tau )\{ T(\partial _{y} ,n)U(y)\} -\{ T(\partial _{y} ,n)\Pi (y,x,\tau )\} ^{*} (U(y)\left. )\right]ds_{y} \right|\le \]

\[\le \int _{S} \left(|\Pi (y,x,\tau )|+|T(\partial _{y} ,n)\Pi (y,x,\tau )|\right)\, \left(|U(y)|+|T(\partial _{y} ,n)U(y)|\right)ds_{y} \le \]

\[\le \delta \int _{S} \left(|\Pi (y,x,\tau )|+|T(\partial _{y} ,n)\Pi (y,x,\tau )|\right)ds_{y} \le C_{3} (x)\delta \tau ^{m} \exp (\tau \, x_{n}^{0} -\tau \, x_{n} ).\]
Then

\[|U(x)|\le C_{4} (x)\tau ^{m} \exp (-\tau \, x_{n} )(M+\delta \exp \tau \, x_{n}^{0} ).\]
Next if we take $\tau =\frac{1}{x_{n}^{0} } ln\frac{M}{\delta } ,$  then we obtain stability estimate:

\[|U(x)|\le C_{4} (x)\delta ^{\frac{x_{n} }{x_{n}^{0} } } \left(ln(\frac{M}{\delta } )\right)^{m} \, .\]

\noindent The theorem is thereby proved.

From proved above theorems we obtain

\textbf{Corollary 1.} \textit{ The limit relation }

\[\mathop{\lim }\limits_{\tau \to \infty } U_{\tau } (x)=U(x),{\kern 1pt} \, \, \, \, \, \, \, \, \, \, \, \, \, \, \, \, \mathop{\lim }\limits_{\delta \to 0} U_{\tau \, \delta } (x)=U(x)\]
hold uniformly on each compact subset of $D.$

 \textbf{3. Regularization of solution of the problem \eqref{GrindEQ__1_}, \eqref{GrindEQ__3_} for a domain of cone type}

Let $x=(x_{1} ,{\rm …},x_{n} )$ and $y=(y_{1} ,{\rm …},y_{n} )$ be points in $E^{n} ,\, \, {\kern 1pt} {\kern 1pt} D_{\rho } $ be a bounded simply connected domain in $E^{n} $ whose boundary consists of a cone surface

\[\Sigma :\quad \alpha _{1} =\tau y_{n} ,{\kern 1pt} \; \alpha _{1}^{2} =y_{1}^{2} +\ldots +y_{n-1}^{2} ,{\kern 1pt} \; \tau =tg\frac{\pi }{2\rho } ,{\kern 1pt} \; y_{n} >0,{\kern 1pt} \; \rho >1\]
and a smooth surface $S,$ lying in the cone. Assume $x_{0} =(0,...0,x_{n} )\in D_{\rho } $.

We constract Karleman matrix. In formula \eqref{GrindEQ__7_},\eqref{GrindEQ__8_} and \eqref{GrindEQ__9_} to take

\[K(\omega )=E_{\rho } \left[\tau \left(\omega -x_{n} \right)\right],\, \, \, \, \, \, \, \tau >0,\, \, \, \, \rho >1.\, \]
Then

\[\Phi (y,x,k)=\Phi _{\tau } (y-x,k),  k>0\]

\[C_{n} \Phi _{\tau } (y-x,k)=\frac{\partial ^{m-1} }{\partial s^{m-1} } \int _{0}^{\infty } Im\, [\frac{E_{\rho } \left(\tau (i\sqrt{u^{2} +s} +y_{n} -x_{n} )\right)}{i\sqrt{u^{2} +s} +y_{n} -x_{n} } ]\frac{\psi (ku)\, du}{\sqrt{u^{2} +s} } \, \, \, \, \, \, (18)\, \]

\[\Phi '_{\tau } (y-x,k)=\frac{\partial \Phi _{\tau } }{\partial \tau } .\]

\[C_{n} \Phi '_{\tau } (y-x,k)=\, \frac{\partial ^{m-1} }{\partial s^{m-1} } \int _{0}^{\infty } Im\left\{\, E'_{\rho } \, \left[\tau (i\sqrt{u^{2} +s} +y_{n} -x_{n} \right]\, \, \right\}\frac{\psi (ku)\, du}{\sqrt{u^{2} +s} } ,\, \, \, \, \, \, \, \, \, \, \, \, \, \, \, (19)\]
where $\quad E_{\rho } (w)-$ Mittag-Löffer`s a entire function [5]. For the functions $\Phi _{\tau } (y-x,k)$ holds Lemma 1 and Lemma 2.

\noindent             Now again to denote by $U_{\tau } (x),{\kern 1pt} \, \, \, \, U_{\tau \, \delta } (x)$ as \eqref{GrindEQ__14_} and \eqref{GrindEQ__17_}. Then holds analogical theorem as \textbf{Theorem} \textbf{1,2,3.}

For    $n=3$ we reduce entirely.

Suppose that $D_{\rho } $ is bounded simple connected domain from $E^{3} $ with boundary consisting of part $\sum  $ of the surface of the cone

\[y_{1}^{2} +y_{2}^{2} =\tau \, y_{3}^{2} ,\quad \tau =tg\frac{\pi }{2\rho } ,\quad \rho >1,\quad y_{3} >0,\]
and of a smooth portion of the surface $S$ lying inside the cone. Assume $x_{0} =(0,0,x_{3} )\in D_{\rho } $.

We construct Carleman`s matrix. In formula \eqref{GrindEQ__7_},\eqref{GrindEQ__8_} we take

\[\Phi _{\tau } (y,x,k)=\frac{1}{4\pi ^{2} E_{\rho } (\tau ^{\frac{1}{\rho } } x_{3} )} \int _{0}^{\infty } Im\frac{E_{\rho } (\tau ^{\frac{1}{\rho } } w)}{i\sqrt{u^{2} +s} +y_{3} -x_{3} } \frac{\cos ku\, du}{\sqrt{u^{2} +s} } ,\, \, \, \, \, \, \, \, \, \, \, \, \, \, \, \, \, (20)\]

where $w=i\sqrt{u^{2} +s} +y_{3} ,\quad E_{\rho } (w)-$ Mittag-Löffer`s a entire function. For the functions $\Phi _{\tau } (y,x,k)$ holds Lemma 1.

If follows from the properties of $E_{\rho } (w)$ that for $y\in \Sigma ,\quad 0<u<\infty $ the function $\Phi _{\tau } (y,x,k)$ defined by \eqref{GrindEQ__18_} its gradient and second partial derivatives

\[\frac{\partial ^{2} \Phi _{\tau } (y,x,k)}{\partial y_{k} \partial y_{j} } ,\quad k,j=1,2,3,\]

tend to zero as $\tau \to \infty $ for a fixed $x\in D_{\rho } .$

Then from \eqref{GrindEQ__7_} we find that the matrix $\Pi (y,x,\tau )$ and its stresses $T(\partial _{y} ,n)\Pi (y,x,\tau )$ also tend to zero as $\tau \to \infty $ on $y\in \Sigma ,$ i.e., $\Pi (y,x,\tau )-$ is the Carleman matrix for the domain $D_{\rho } $ and the part $\Sigma $ of the boundary.

For the $U(x)-$ regular solution system \eqref{GrindEQ__1_} following integral formula holds

\[U(x)=\int _{\partial D_{\rho } } [\Pi (y,x,\tau )\{ T(\partial _{y} ,n)U(y)\} -\{ T(\partial _{y} ,n)\Pi (y,x,\tau )\} ^{*} U(y)]ds_{y} .(21)\]

By $x\in D_{\rho } $ we denote $U_{\tau } (x)$ follows:

\[U_{\tau } (x)=\int _{S} [\Pi (y,x,\tau )\{ T(\partial _{y} ,n)U(y)\} -\{ T(\partial _{y} ,n)\Pi (y,x,\tau )\} ^{*} U(y)]ds_{y} .(22)\]

The following theorem holds.

\textbf{Theorem 4.} \textit{ Let $U(x)$ be a regular solution of system \eqref{GrindEQ__1_} in $D_{\rho } $ such that }

\[|U(y)|+|T(\partial _{y} ,n)U(y)|\le M,{\kern 1pt} \; y\in \Sigma .\, \, \, \, \, \, \, \, \, \, \, \, \, \, \, \, \, \, \, \, \, \, \, \, \, \, \, \, \, \, \, \, \, \, \, \, \, \, \, \, \, \, \, \, \, \, \, (23)\]

Then for $\tau \ge 1$ the following estimate is valid:

\[|U(x_{0} )-U_{\tau } (x_{0} )|\le MC_{\rho } (x_{0} )\tau ^{3} \exp (-\tau \, x_{3}^{\rho } ),\]

where $x_{0} =(0,0,x_{3} )\in D_{\rho } ,\quad x_{3} >0,\quad $

\[C_{\rho } (x_{0} )=C_{\rho } \int _{\Sigma } \frac{1}{r_{0}^{3} } ds_{y} ,\quad r_{0} =|y-x_{0} |,\quad C_{\rho } -constant.\]

\textbf{Proof.} By analogy with proved Theorem 2 and Theorem 3 from \eqref{GrindEQ__18_} and \eqref{GrindEQ__19_} we obtain

\[|U(x_{0} )-U_{\tau } (x_{0} {\rm \; \; }|\le M\int _{\Sigma } [|\Pi (y,x_{0} ,\tau )|+|T(\partial _{y} ,n)\Pi (y,x_{0} ,\tau )|]ds_{y} .\, \, \, \, \, \, \, \, \, \, \, \, \, \]

By formula \eqref{GrindEQ__20_} we have following inequality:

\[\left|\Phi _{\tau } (y,x,k)\right|\le C_{\rho }^{(1)} E_{\rho }^{-1} (\tau ^{\frac{1}{\rho } } \gamma )r^{-1} ,\]

\[\left|\frac{\partial \Phi _{\tau } (y,x,k)}{\partial y_{i} } \right|\le C_{\rho }^{(2)} \tau E_{\rho }^{-1} (\tau ^{\frac{1}{\rho } } \gamma )r^{-2} \]

\[\left|\frac{\partial ^{2} \Phi _{\tau } (y,x,k)}{\partial y_{k} \partial y_{j} } \right|\le C_{\rho }^{(3)} \tau ^{2} E_{\rho }^{-1} (\tau ^{\frac{1}{\rho } } \gamma )r^{-3} .\]

Then from \eqref{GrindEQ__7_}

\[|\Pi (y,x,\tau )|\le C_{\rho }^{(4)} \tau ^{2} E^{-1} (\tau ^{\frac{1}{\rho } } \gamma )r^{-3} ,\]

\[|T(\partial _{y} ,n)\Pi (y,x,\tau )|\le C_{\rho }^{(5)} \tau ^{3} E^{-1} (\tau ^{\frac{1}{\rho } } \gamma )r^{-4} .\]

Therefore  we obtain

\[|U(x_{0} )-U_{\tau } (x_{0} )|\le MC_{\rho } (x_{0} )\tau ^{3} \exp (-\tau \, x_{3}^{\rho } ),\]

where

\[C_{\rho } (x_{0} )=C_{\rho } \int _{\Sigma } \frac{1}{r_{0}^{3} } ds_{y} ,\quad r_{0} =|y-x_{0} |,\quad C_{\rho } -constant.\]

The theorem is thereby proved.

Suppose that instead of $U(y)$ and $T(\partial _{y} ,n)U(y)$ gives their continuous approximations $f_{\delta } (y)$ and $g_{\delta } (y)$ such that

\[\mathop{\max }\limits_{S} |U(y)-f_{\delta } (y)|+\mathop{\max }\limits_{S} |T(\partial _{y} ,n)U(y)-g_{\delta } (y)|\le \delta ,{\kern 1pt} \; \, \, \, 0<\delta <1.\]

Define the function $U_{\tau \, \delta } (x)$ by

\[U_{\tau \, \delta } (x)=\int _{S} [\Pi (y,x,\tau )g_{\delta } (y)-\{ T(\partial _{y} ,n)\Pi (y,x,\tau )\} ^{*} f_{\delta } (y)]ds_{y} ,\, \, \, \, \, \, \, \, \, \, \, \, \, \]

The following theorem holds

\textbf{Theorem 5.} \textit{ Let $U(x)$ is a regular solution of system \eqref{GrindEQ__1_} in the domain $D_{\rho } $ satisfying the condition \eqref{GrindEQ__23_}, then }

\[|U(x_{0} )-U_{\tau \, \delta } (x_{0} )|\le C_{\rho } (x_{0} )\delta ^{q} (ln\frac{M}{\delta } )^{3} ,\]

where $\tau =(\tau R)^{-\rho } ln\frac{M}{\delta } ,\quad R^{\rho } =\mathop{\max }\limits_{S} Re(i\sqrt{s} +y_{3} )^{\rho } ,$

\[q=(\frac{\gamma }{R} )^{\rho } ,\quad C_{\rho } (x_{0} )=C_{\rho } \int _{\Sigma } \left[\frac{1}{r_{0}^{3} } +\frac{1}{r_{0}^{4} } \right]ds_{y} .\]

 The proof theorem is similar to those of Theorem 3 and 4.

\textbf{Corollary 2.} \textit{ The limit relation }

\[\mathop{\lim }\limits_{\tau \to \infty } U_{\tau } (x)=U(x),{\kern 1pt} \, \, \, \, \, \, \, \, \, \, \, \, \, \, \, \, \mathop{\lim }\limits_{\delta \to 0} U_{\tau \, \delta } (x)=U(x)\]

hold uniformly on each compact subset of $D_{\rho } .$

 \textbf{References}

[1]. M.M.Lavrent'ev. Some Ill-Posed Problems of Mathematical Physics [in Russian], Computer Center of the Siberian Division of the Russian Academy of Sciences, Novosibirck (1962) 92p.

[2]. I.G.Petrovskii. Lectures on Partial Differential Equations [in Russians],Fizmatgiz,Moscow,(1961).

[3]. V.D.Kupradze, T.V.Burchuladze, T.G.Gegeliya, ot.ab. Three-Dimensional Problems of the Mathematical Theory of Elasticity and ... [in Russian],Nauka,Moscow,1976.

[4]. Sh.Ya.Yarmukhamedov. Dokl.Acad.Nauk SSSR [Soviet Math.Dokl.], V.357, No.3, p.320-323.(1997).

 [5]. M.M. Dzharbashyan, Integral Transformations and Representations of Functions in a Complex Domain [in Russian],Nauka, Moscow.1966.

\noindent

\end{document}